\documentclass[journal]{IEEEtran}

\usepackage{graphicx}
\usepackage{cite}
\usepackage{amsmath, bm, graphicx, xcolor, makecell}

\allowdisplaybreaks

\ifCLASSOPTIONcompsoc
	\usepackage[caption=false,font=normalsize,labelfont=sf,textfont=sf]{subfig}
\else
	\usepackage[caption=false,font=footnotesize]{subfig}
\fi
\graphicspath{{images/}}

\begin{document}

\title{A Trilevel Model for Segmentation of the Power Transmission Grid Cyber Network}

\author{Bryan~Arguello,
        Emma~S.~Johnson,
        and~Jared~L.~Gearhart%
\thanks{
        This work was supported by Sandia National Laboratories' Laboratory Directed Research and Development (LDRD) program.
        Sandia National Laboratories is a multimission laboratory
        managed and operated by National Technology and Engineering
        Solutions of Sandia, LLC, a wholly owned subsidiary
        of Honeywell International, Inc., for the U.S. Department
        of Energy’s National Nuclear Security Administration under
        contract DE-NA0003525. SAND NO. 2021-10208 O.
        The views expressed in the article do not necessarily represent the views of the U.S. Department of Energy or the United States Government.}%
\thanks{B. Arguello and J. Gearhart are with Sandia National Laboratories, Albuquerque, NM 87185 USA (email: \{barguel; jlgearh\}@sandia.gov)}%
\thanks{E. Johnson is with the H. Milton Stewart School of Industrial and Systems Engineering, Georgia Institute of Technology, Atlanta, GA 30332 USA, and also with Sandia National Laboratories, Albuquerque, NM 87185 USA (email: ejohnson335@gatech.edu)}%
}

{}

\maketitle

\begin{abstract}
Network segmentation of a power grid's communication system is one way to make the grid more resilient to cyber attacks. 
We develop a novel trilevel programming model to optimally segment a grid communication system, taking into account the actions of an information technolology (IT) administrator, attacker, and grid operator. 
The IT administrator is given an allowance to segment existing networks, and the attacker is given a fixed budget to attack the segmented communication system in an attempt to inflict damage on the grid. Finally, the grid operator is allowed to redispatch the grid after the attack in order to minimize damage. 
The resulting problem is a trilevel interdiction problem, which we solve by leveraging current research in bilevel branch and bound. 
We demonstrate the benefits of optimal network segmentation through case studies on the 9-bus WSCC system and the 30-bus IEEE system. 
These examples illustrate that network segmentation can significantly reduce the threat posed by a cyber attacker with perfect knowledge of the grid.
\end{abstract}

\begin{IEEEkeywords}
Power system protection, Cyberspace, Networked control systems, Optimization, Multilevel systems
\end{IEEEkeywords}

\section*{Nomenclature}

\subsection*{Sets}

\begin{description}
	\item [$\mathcal E$] Security enclaves 
	\item [$\mathcal B$] Balancing authority entities
	\item [$\mathcal C$] Control center entities
	\item [$\mathcal S$] Substation entities
	\item [$\mathcal T$] Tiers of entities in the communication network: \newline $\mathcal T = \{\mathcal B, \mathcal C, 
	\mathcal S\}$
	\item [$\mathcal Z$] Pairs of adjacent tiers in communication network: $\mathcal Z = \{(\mathcal B, \mathcal C), (\mathcal C, \mathcal S)\}$. Generic members of $\mathcal Z$ are always denoted $(A, B)$ where $A$ is the tier above and $B$ is the tier below.
	\item [$\mathcal E_0(T)$] Existing enclaves in tier $T \in \mathcal T$
	\item [$\mathcal E_1(T)$] New enclaves in tier $T \in \mathcal T$
	\item [$\mathcal E(T)$] Enclaves in tier $T \in \mathcal T$: $\mathcal E(T) = \mathcal E_0(T) \cup \mathcal E_1(T)$
	\item [$T_n$] Entities in tier $T$ which communicate with entity $n$ one tier above
	\item [$\mathcal{K}$] Transmission lines
	\item [$\mathcal{G}$] Generators
	\item [$\mathcal{G}_s$] Generators at substation $s$
	\item [$\mathcal L$] Loads
	\item[$\mathcal{L}_s$] Loads at substation $s$
	\item [$\mathcal R$] Relays
	\item [$\mathcal R_s$] Relays at substation $s$
	\item [$\mathcal R_d$] Relays which control load $d$
	\item [$\mathcal R_k$] Relays which control line $k$
	\item [$\mathcal R_g$] Relays which control generator $g$
\end{description}

\subsection*{Parameters}
\begin{description}
	\item[$Q_{n,e}$] Binary indicating if existing enclave $e$ is in entity $n$
	\item[$U$] Maximum number of enclaves attacker can penetrate
	\item[$\Theta_k$] Transformer shift angle on line $k$
	\item[$B_k$] Line charging susceptance of line $k$
	\item[$D_d$] Demand at load $d$
	\item[$s(d)$] Substation served by load $d$
	\item[$s(g)$] Substation served by generator $g$
	\item[$\overline P_g$] Maximum real power output of generator $g$
	\item[$\overline F_k$] Thermal limit of line $k$
	\item[$o(k)$] Origin bus of line $k$
	\item[$d(k)$] Destination bus of line $k$
\end{description}

\subsection*{Binary Variables}
IT Operator Decisions:
\begin{description}
	\item[$x_{e,r}$] Indicates whether or not enclave $e \in \mathcal E(\mathcal S)$ communicates with relay $r$
	\item[$y_{e,f}$] Indicates whether or not enclave $e$ communicates with enclave $f$
	\item[$q_{n,e}$] Indicates whether or not enclave $e$ is in entity $n$
	\item[$t_{e,n}$] Indicates whether or not enclave $e$ communicates with entity $n$ (in the tier below)
\end{description}

Attacker decisions:
\begin{description}
	\item[$z_e$] Indicates whether or not the attacker infiltrates enclave $e$
	\item[$\delta_r$] Indicates whether or not relay $r$ is compromised
	\item[$v_k$] Indicates whether or not line $k$ is operational
	\item[$u_d$] Indicates whether or not the load $d$ is online
	\item[$w_g$] Indicates whether or not generator $g$ is operational
\end{description}

\subsection*{Continuous Variables}
Grid operator decisions:
\begin{description}
	\item[$\theta_s$] Voltage angle at substation $s$
	\item[$p_g$] Real power output of generator $g$
	\item[$f_k$] Real power flow on line $k$
	\item[$l_d$] Real power load shed at load $d$
	\item[$L$]   Real power total load shed
\end{description}

\pagebreak

\section{Introduction}

\IEEEPARstart{T}{he} electric grid, like many infrastructure systems, historically relied on air-gapped and specialized control networks for security. For convenience and cost savings, grid/utility owners are increasingly integrating traditional information technologies (IT) into their control systems. 
This shift has increased the attack surface for electric grids and has made them more vulnerable to cyber attacks. This was illustrated in the December 2015 cyber attack against the Ukrainian power grid \cite{cisa2016}. 
According to \cite{dni2019}, nation-state actors in particular pose a significant threat to cyber-physical systems and also possess the ability to gather the necessary information in advance of attacks. 
Criminal organizations, terrorists, hackers, and hackivists may also pose a threat to these systems \cite{gao2019}.

Network segmentation is one strategy that has been proposed to improve cyber network security \cite{nsa2016}. Dividing networks into small segments and restricting communications between segments can limit the scope of attacks and the attackers' ability to pivot within a network.  
Several governmental and regulatory agencies have provided guidance that motivates segmentation and offer strategies for implementing it \cite{dhs2016}, \cite{pillitteri2014}, \cite{stouffer2015}, \cite{acsc2020}. 
However, as noted in \cite{wagner2016}, this guidance tends to provide general design principles and architectures, but not specific recommendations or tools for cyber network designers. 

Designing segmentation strategies for cyber-physical systems presents a particularly difficult challenge as it requires consideration of two complex interconnected systems. 
In an IT setting, segmentation decisions can be made based on business functions and the expected consequences of disruptions to these functions (e.g. interruption to payment systems, loss of customer data). 
In a grid setting, the cyber and physical layers are connected and distributed over wide areas; coordinated attacks create the potential for consequences that are greater than the sum of their parts. 
A key recommendation in \cite{gao2019} is the need to ``evaluate the potential risk of a coordinated cyberattack on geographically distributed targets."
 
In this paper, we explore the use of optimization to identify strategies for segmenting cyber networks. 
We use an explicit representation of the cyber and grid network layers to capture the interconnected nature of these systems and the resulting impacts of attacks. 
We also employ trilevel optimization techniques to capture the adaptive nature of and interactions between cyber defenders, attackers, and infrastructure operators. 
Finally, given the potential capabilities of the attacker, we assume that the attacker has full knowledge of the system when crafting attacks. 

\subsection{Literature Review}
Several researchers have previously considered applying optimization to cyber security problems that are related to network segmentation.
In\cite{jeong2005}, a heuristic optimization algorithm is used to place the minimum number of intrusion detection systems (IDS) on the nodes of a cyber network while ensuring that no nodes are more than $n$ hops from the nearest detection system. 
This heuristic approach is also applied by \cite{zhang2013} to place trust systems on the smart grid cyber network.

In \cite{gonzalez2011}, trust system placement is combined with network segmentation. 
A mixed-integer linear program (MILP) is used to partition a power system supervisory control and data acquisition (SCADA) network into a collection of enclaves such that a trust node is placed on at least one end of all arcs connecting different enclaves. 
This ensures that all communications between segments are scanned. 
The authors also develop a heuristic algorithm for partitioning and placing trust nodes on larger networks. 

In \cite{wagner2016}, simulated annealing is used to identify segmentation strategies for a computer network in an air operations center that supports 40 military missions. 
An emulation test bed is used to simulate the performance of each segmentation strategy against a ``dumb" worm attack that propagates in a manner similar to a disease in a Susceptible-Exposed-Infectious-Recovered (SEIR) model. The quality of each segmentation strategy is based on mission delay and the availability of devices on the network. 
In \cite{hemberg2018}, this approach is extended to an algorithm that co-evolves the parameters available to the attacker and defender in order to optimize both of their strategies. 
The attacker still employs a ``dumb" worm attack and cannot control how it propagates, but they tune the intensity of the attack to balance the number of components compromised against the probability of detection.

In  \cite{hasan2016a}, a heuristic method that uses minimum spanning trees is used to improve performance and scalability for the trust node placement problem on smart grid SCADA systems. 
Variants of this model are considered in \cite{hasan2016b} and \cite{hasan2016c}. 
These variants account for node centrality and communication latency when deciding where to place trust nodes. 
In \cite{hasan2019}, this approach is further expanded to consider both link coverage (the number of links adjacent to a trust node) and path tolerance (the longest path not passing through a trust node) when making placement decisions. 

In all of the previously mentioned references that focus on cyber-physical security, only the cyber network is considered, while effects on the underlying physical system it controls are not modeled. 
There are several examples where segmentation decisions are made based on the underlying physical system managed by the cyber network. 
In \cite{genge2012}, an emulation test bed is used to compare two segmentation strategies developed by subject matter experts (SME's) for SCADA systems in a chemical plant. 
They demonstrate that in process-flow operations, security and safety gains can be made by separating SCADA devices for in-flow and out-flow values. 
This approach is extended in \cite{genge2014} to include an SME-informed heuristic approach that creates segmentation strategies to separate control of in-flows and out-flows. 
In \cite{arief2020}, network segmentation is used to prevent domino effects in storage tank facilities (e.g. chained explosions of nearby tanks) by an attacker who can attack at most one segment. 
Belief networks are used to estimate the probability of chained reactions and a graph centrality measure is used to find the segmentation strategy.
These works use testbeds to study the effects of network segmentation on the underlying physical system. 
However, they do not use optimization to make their network segmentation decisions.

Prior research has used bilevel programming to consider the effects on the power grid from an intentional attack, though these model only the physical grid, not the cyber communication system.
For example, \cite{salmeron2009worst} and \cite{motto2005mixed} formulate and solve a bilevel program that models an attacker with perfect knowledge of the grid and the response of a grid operator. 
The attacker directly de-energizes power grid components such as transmission lines and generators. 
The grid operator then responds through generator redispatch and load shed to minimize unmet demand. 
This attacker-defender model is adapted for a cyber attacker in \cite{CastilloACS2019}, though still without modeling the cyber network explicitly.

Finally, several authors have employed trilevel programming to preemptively protect the power grid against malicious attacks in \cite{delgadillo2009analysis}, \cite{wu2016efficient}, \cite{alguacil2014trilevel}, and \cite{YuanZZ2014}. 
These works focus on hardening power grid components such as transmission lines, generators, and buses against optimal physical attacks on the grid. 
However, they do not consider the grid's communication system in modeling these attacks.

\subsection{Contributions and Paper Organization}
In this work, we propose segmenting the grid's communication network optimally by explicitly modeling both the cyber network and the physical network.
Specifically, we use trilevel optimization to model the interconnected decisions made by an IT administrator segmenting a grid's communication network, a cyber attacker maliciously trying to damage the grid by attacking the segmented network, and a grid operator mitigating damage from the attack. 
The main contributions of this paper are:
\begin{itemize}
	\item The formulation of a trilevel programming model to identify provably optimal segmentation strategies for a power transmission cyber-physical system which reduce the severity of the worst attack available to a malicious cyber attacker.
	\item The use of an infrastructure model (i.e., a DC optimal power flow (DCOPF) grid model) to inform network segmentation decisions.
	\item A description of how to reformulate the trilevel model as a bilevel program, which can be solved using algorithms from the literature (e.g. bilevel branch-and-bound).
	\item A case study showing network segmentation results for both the 9-bus WSCC and 30-bus IEEE test cases, where segmentation reduces the severity of the worst-case attack by 31\% and 56\%, respectively.
\end{itemize}
In the remainder of the paper, Sections \ref{sec:model-high-level} and \ref{sec:model-formulation} give a description of the model. 
Section \ref{sec:methodology} details our solution methodology. 
Section \ref{sec:results} presents our case studies on the 9-bus and 30-bus IEEE systems. 
We conclude in Sections \ref{sec:future-work} and \ref{sec:conclusion} with ideas for future work and conclusions.

\section{Cyber-Physical Network Segmentation Model}\label{sec:model-high-level}
We first describe the trilevel network segmentation model at a high level.
We give the mathematical formulation in the next section. 
When we talk about the cyber-physical system, we draw a distinction between {\it entities}, the bodies which control the grid, {\it enclaves}, the networks that make up an entity's IT infrastructure, and {\it physical components} such as the lines, generators, and loads on the physical grid.
While our model is more general, in this paper we consider three types of entities: balancing authorities, control centers, and substations.

\subsection{Modeling the Communication Network}
For a description of how a power transmission grid communication system network is structured, we refer the reader to \cite{gaudet2020firewall}.
We simplify the communication network to a 3-tier forest where each substation enclave, control center enclave, and balancing authority enclave is represented as a node, and parent-child relationships represent that the parent enclave communicates with the child enclave (i.e., data flows between the parent and child enclave).
Every substation enclave is the child of exactly one control center enclave. 
In turn, every control center enclave is the child of exactly one balancing authority enclave. 
\begin{figure*}
	\centering	
	\subfloat[The 9-bus system communication network before segmentation: Arcs indicate that the parent component communicates with the child component. 
	The balancing authority, control center, and substation enclaves comprise the communication system, and its control over the relays is shown as the leaves of the tree. 
	(Note that in this example, we include fewer relays than in the 9-bus system used in the case study in Section \ref{subsec:9-bus-results}.) ]{\includegraphics[width=0.48\linewidth]{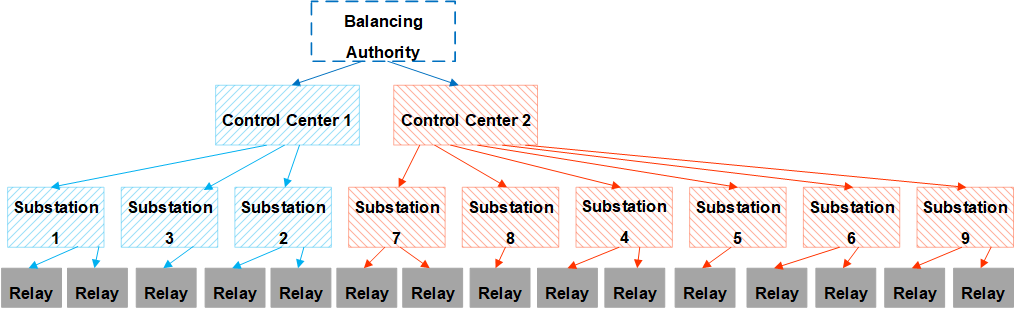}\label{fig:9bus-graph-before-segmentation}} \hspace{0.25cm}
	\subfloat[The 9-bus system communication network after segmentation. 
	The balancing authority enclave has been subdivided into two enclaves, each of which communicates with one of the control centers. 
	The contol centers also have two enclaves each, as does Substation 4.]{\includegraphics[width=0.48\linewidth]{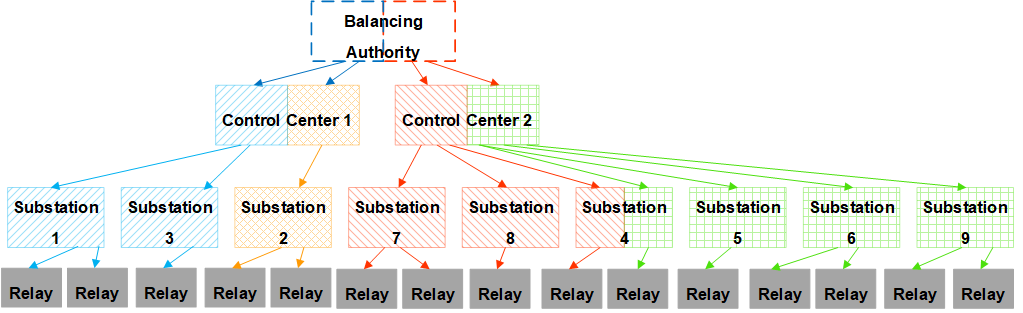}\label{fig:9bus-graph-after-segmentation}} \\
	\subfloat[The 9-bus physical network before segmentation: Grid components are shaded based on what control-center enclave in the communication network controls them, and the dotted line represents what balancing authority enclave controls the components it contains. 
	Before segmentation, there is only one balancing authority enclave. 
	Buses 1, 2, and 3 are all controlled by relays controlled by Substations 1, 2, or 3, and hence controlled by Control Center 1. The rest of the network is controlled by Control Center 2. ]{\includegraphics[width=0.48\linewidth]{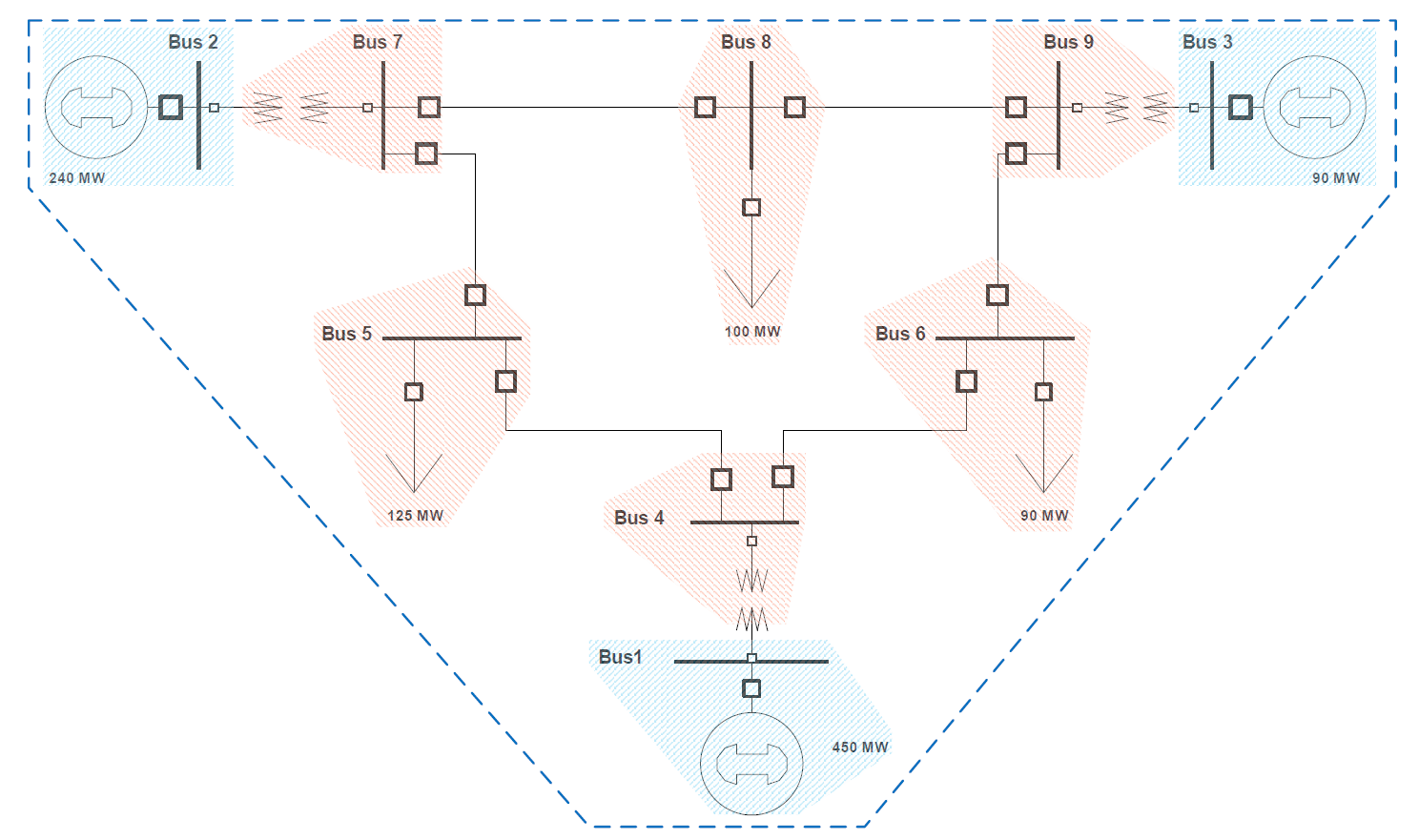}\label{fig:9bus-network-before-segmentation-noattack}} \hspace{0.25cm}
	\subfloat[The 9-bus physical network after segmentation: Now the network is divided between two balancing authority enclaves, and each of the two sets of components controlled by Control Centers 1 and 2 are now divided between two different enclaves in their control center. 
	Dotted angled lines inside of the shaded regions indicate grid components which are on different relays. 
	For example, two of the lines adjacent to Bus 4 are on a relay controlled by the green (cross-hatched) enclave of Control Center 2, and the other line is controlled by the red (striped) enclave. ]{\includegraphics[width=0.48\linewidth]{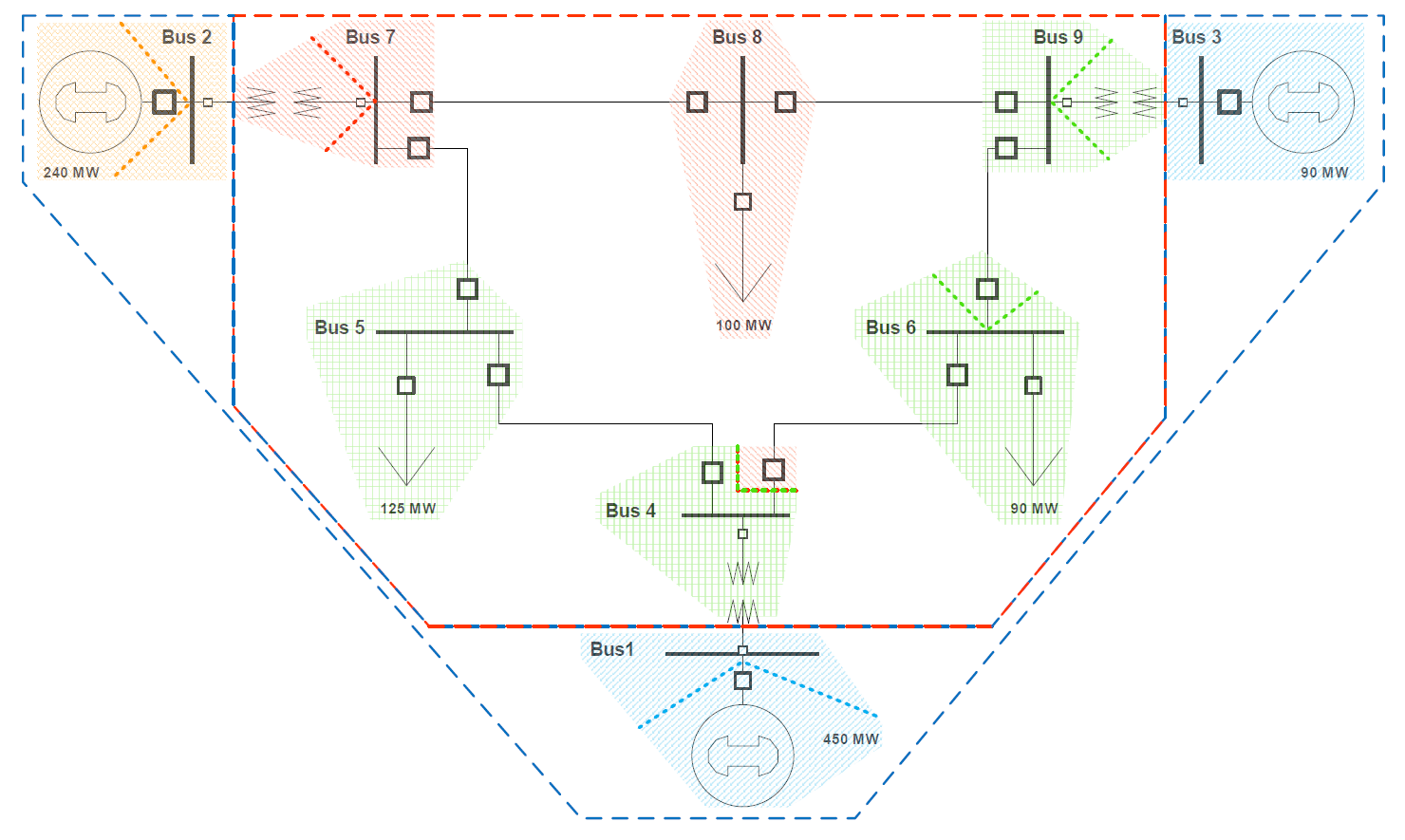}\label{fig:9bus-network-after-segmentation-noattack}}
	\caption{The 9-bus system and a graph representation of its communication network, both before and after segmentation.
		The communication network is shown before segmentation in (a) and after segmentation in (b).
		The physical grid is shown before segmentation in (c) and after segmentation in (d).
		Colors are used to map which physical components are controlled by which control-center enclaves in the communication network.}
	\label{fig:9bus-noattack}
\end{figure*}
We show an example of the forest representation of the 9-bus WSCC test system communication network in Fig.~\ref{fig:9bus-graph-before-segmentation}.
In this example, every entity has only one enclave. 
The balancing authority enclave communicates with two different control center enclaves.
The first of these communicates with three substation enclaves, and the second with six substation enclaves. 

The physical power grid is connected to its communication system through relays, which are drawn as the leaves of the tree in Fig.~\ref{fig:9bus-graph-before-segmentation}. 
Individual relays are controlled by exactly one substation enclave and can be used by an attacker to de-energize transmission lines, generators, and loads. 
In Fig.~\ref{fig:9bus-network-before-segmentation-noattack}, each shaded block of grid components is colored to indicate which control center enclave from Fig.~\ref{fig:9bus-graph-before-segmentation} ultimately controls it. 

\subsection{Modeling Network Segmentation}\label{subsec:modeling-segmentation}
In the outermost problem of the trilevel model, we model the decisions of an IT administrator segmenting the network in order to minimize the load shed from a worst-case cyber attack.
We model network segmentation by allowing the IT administrator to partition any security enclave into two or more security enclaves while still respecting the underlying control structure among the entities. 
By only subdividing existing assignments, we aim to generate new segmentation strategies that are compatible with the original topology. 
In the graph representation, segmentation corresponds to expanding the graph by dividing nodes.
If a node is divided, the new subnetwork must be assigned to the same entity as the original node.
For example, if a substation enclave is divided, both new enclaves are subnetworks of the same substation entity.
Edges are added to the expanded graph according to the following rules:
\begin{enumerate}
	\item A new enclave may only communicate with (i.e., be a child of) a parent enclave if the entities containing the new enclave and the parent enclave communicated in the original graph. As an example, a substation enclave may only communicate with a control center enclave if the corresponding substation and control center originally communicated before segmentation.
	\item After partitioning a substation enclave, that substation's relays must be reassigned so that each relay is controlled by exactly one of the new substation enclaves. The result of this partitioning is that a substation's grid connections are separated by the substation's enclaves, making that substation less vulnerable overall.
	\item Every substation enclave must communicate with at least one relay.
\end{enumerate}
Continuing our example, see Fig.~\ref{fig:9bus-graph-after-segmentation} for a possible segmentation of the 9-bus network. 
We see that, after network segmentation, the balancing authority has two enclaves, one for each control center.
Both control centers gain a new enclave. For each control center, the substations it originally controlled are divided between its two new enclaves.
Finally, Substation 4 is segmented into two enclaves, each communicating with a separate enclave within Control Center 2. 
Note that each of the two relays controlled by Substation 4 is assigned to a different enclave. 
The corresponding diagram of the physical grid is shown in Fig.~\ref{fig:9bus-network-after-segmentation-noattack}.

\subsection{Attacker-Defender Model}
In the inner two problems of the trilevel model, we model a cyber attacker who finds the highest-load-shed attack on the segmented network given that the grid operator can redispatch after their attack in order to minimize load shed. 
We assume that attacks begin by infiltrating balancing authority enclaves, then progress by gaining access to the control center enclaves followed by substation enclaves, until the attack reaches the relays, where it can produce an effect on the grid.
At each step, the attacker can never infiltrate a node in the graph without first infiltrating its parent.
We assume the attacker has an upper limit on the number of enclaves that can be compromised, and every time the attacker gains access to a security enclave, exactly one unit of this budget is consumed. 
In other words, one unit of budget is used for every node in the graph the attacker reaches.
When the attacker gains access to a substation enclave, we assume that the attacker will use all relays controlled by that substation enclave to de-energize every accessible grid component.
To gain some intuition about the benefits of network segmentation, note for example that originally in Fig.~\ref{fig:9bus-noattack}, an attacker would need to access three enclaves to gain control of the two relays controlled by Substation 4. 
After segmentation, the attacker would need to access five enclaves to gain control of the same relays.

Finally, given an attack, the grid operator solves a DCOPF on the remaining components to minimize total load shed. 

\section{Trilevel Formulation}\label{sec:model-formulation}

We now give a mathematical formulation of the problem described in Section \ref{sec:model-high-level}. 
Because many constraints involving (parent, child) entity pairs hold at multiple levels -- both (balancing authority, control center) and (control center, substation) -- we use the set $\mathcal T$ to denote the set of all three communication network entities and the set $\mathcal Z = \{(\mathcal B, \mathcal C), (\mathcal C, \mathcal S)\}$ to denote adjacent tiers of entities.
The trilevel model is:
\begin{align}\label{trilevel-model}
\min_{(x, y, q, t) \in \mathcal D} \quad \max_{(\delta, z, u, v, w) \in \mathcal A(x,y)} \quad  \min_{(\theta, f, p, l) \in \mathcal O(u, v, w)} \sum_{d \in \mathcal L} l_d,
\end{align}
where $\mathcal D$ is the feasible region of the network designer, $\mathcal A(x,y)$ is the feasible region of the attacker, given the design decisions $x$ and $y$, and $\mathcal O(u, v, w)$ is the DCOPF feasible region given the attack defined by $u$, $v$, and $w$.

We model $\mathcal D$ with the following constraints:
\begin{align}
&\sum_{r \in \mathcal R} x_{e,r} \geq 1, \qquad \forall e \in \mathcal E(\mathcal S) \label{substation-subnets-talk-to-a-relay}\\
&\sum_{e \in \mathcal E(\mathcal S)} x_{e,r} = 1, \qquad \forall r \in \mathcal R \label{each-relay-talks-to-a-subnet} \\
& q_{s,e} \leq \sum_{r \in \mathcal R_s} x_{e,r}, \qquad \forall s \in \mathcal S, e \in \mathcal E_1(\mathcal S) \label{new-substation-subnet-talks-to-some-relay} \\
& q_{s,e} \geq x_{e,r}, \qquad \forall s \in \mathcal S,  r \in \mathcal R_s, e \in \mathcal E_1(\mathcal S) \label{no-relays-means-no-substation-subnet-new}\\
& Q_{s,e} \leq \sum_{r \in \mathcal R_s} x_{e,r}, \qquad \forall s \in \mathcal S, e \in \mathcal E_0(\mathcal S) \label{existing-substation-subnet-talks-to-some-relay} \\
& Q_{s,e} \geq x_{e,r}, \qquad \forall s \in \mathcal S,  r \in \mathcal R_s, e \in \mathcal E_0(\mathcal S) \label{no-relays-means-no-substation-subnet-existing}\\
& \sum_{e \in \mathcal E(A)} y_{e,f} = 1, \qquad \forall (A, B) \in \mathcal Z, f \in \mathcal E(B) \label{all-subnets-mapped} \\
&\sum_{n \in T} q_{n,e} = 1, \qquad \forall T \in \mathcal T, e \in \mathcal E_1(T) \label{each-subnet-assigned-one-entity} \\
&\begin{aligned}
t_{e,n} \leq \sum_{f \in \mathcal E_0(B)} y_{e,f}Q_{n,f} + \sum_{f \in \mathcal E_1(B)} y_{e,f} q_{n,f}, \\ 
\forall (A, B) \in \mathcal Z, e \in \mathcal E(A), n \in B
\end{aligned} \label{inter-layer-communication-1} \\
&\begin{aligned}
t_{e,n} \geq y_{e,f}Q_{n,f}, \qquad \forall &(A, B) \in \mathcal Z, e \in \mathcal E(A), \\
&n \in B, f \in \mathcal E_0(B)
\end{aligned} \label{inter-layer-communication-2} \\
&\begin{aligned}
t_{e,n} \geq y_{e,f}q_{n,f}, \qquad \forall &(A, B) \in \mathcal Z, e \in \mathcal E(A), \\
&n \in B, f \in \mathcal E_1(B)
\end{aligned} \label{inter-layer-communication-3} \\
& \begin{aligned}
t_{e,n} \leq Q_{m,e}, \qquad \forall &(A, B) \in \mathcal Z, e \in \mathcal E_0(A), \\
&m \in A, n \in B_m
\end{aligned} \label{communication-implies-belonging-existing}\\
& \begin{aligned}
t_{e,n} \leq q_{m,e}, \qquad \forall &(A, B) \in \mathcal Z, e \in \mathcal E_1(A), \\
&m \in A, n \in B_m
\end{aligned} \label{communication-implies-belonging-new}\\
&y_{e,f} \in \{0, 1\}, \forall (e,f) \in (\mathcal E(\mathcal C) \times \mathcal E(\mathcal S)) \cup (\mathcal E(\mathcal B) \times \mathcal E(\mathcal C)) \label{y-domain} \\
&x_{e,r} \in \{0, 1\}, \qquad  e \in \mathcal E(\mathcal S), \forall r \in \mathcal R \label{x-domain} \\
&q_{n,e} \in \{0, 1\}, \qquad \forall (n,e) \in \cup_{T \in \mathcal T} ( T \times \mathcal E_1(T)) \label{q-domain} \\
&t_{e,n} \in \{0,1\}, \quad \forall (e,n) \in (\mathcal E(\mathcal C) \times \mathcal S) \cup (\mathcal E(\mathcal B) \times \mathcal C) \label{t-domain}
\end{align}
In (\ref{substation-subnets-talk-to-a-relay}), we require that every substation enclave controls at least one relay, and in (\ref{each-relay-talks-to-a-subnet}) we enforce that every relay is controlled by exactly one substation enclave.
Referring back to the graph in Fig.~\ref{fig:9bus-noattack}, the $x_{r,e}$ variables represent whether or not there is an arc from enclave $e$ in the substation tier to relay $r$.
Constraints (\ref{new-substation-subnet-talks-to-some-relay})-(\ref{no-relays-means-no-substation-subnet-existing}) enforce that an enclave $e$ is assigned to a substation if and only if some relay controlled by that substation is controlled via $e$.

In (\ref{all-subnets-mapped}), we require that, between any two adjacent tiers, every enclave in the tier below is controlled by exactly one enclave in the tier above.
Since $y_{e,f}$ determines whether or not there is an arc from enclave $e$ to enclave $f$, these constraints enforce that substation enclaves and control center enclaves have exactly one parent.

Constraints (\ref{each-subnet-assigned-one-entity}) guarantee that, for each of the tiers of the communication network, every new enclave is assigned to exactly one entity (substation, control center, or balancing authority depending on the tier). 
Constraints (\ref{inter-layer-communication-1})-(\ref{inter-layer-communication-3}) enforce that enclave $e$ communicates with entity $n$ in the layer below $e$ if and only if entity $n$ has an enclave $f$ that communicates with $e$.
Returning to the graphs in Fig.~\ref{fig:9bus-noattack}, $t_{e,n}$ is an indicator of whether enclave $e$ is a parent of some enclave assigned to entity $n$.
Constraints (\ref{communication-implies-belonging-existing})-(\ref{communication-implies-belonging-new}) enforce that if enclave $e$ communicates with entity $n$ in the network layer below it, then $e$ must be assigned to an entity $m$ which controls $n$. Note that $B_m$ is an instance of the set $T_n$ with $T=B$ and $n=m$.
Finally, (\ref{y-domain})-(\ref{t-domain}) give variable domains.

Note that (\ref{inter-layer-communication-1}) and (\ref{inter-layer-communication-3}) contain bilinear terms. Since all variables in these products are binaries, they can be reformulated linearly with the introduction of one new binary variable per bilinear term.
That is, we introduce $\beta_{e, f, n} \in \{0,1\}$ for all $(A, B) \in \mathcal Z, e \in \mathcal E(A),f \in \mathcal E_1(B), n \in B$ and the constraints
\begin{align}
&\begin{aligned}
\beta_{e, f, n} \leq y_{e,f}, \qquad \forall &(A, B) \in \mathcal Z, e \in \mathcal E(A), \\
&f \in \mathcal E_1(B), n \in B 
\end{aligned}\label{mccormick1}\\
&\begin{aligned}
\beta_{e, f, n} \leq q_{n,f}, \qquad  \forall &(A, B) \in \mathcal Z, e \in \mathcal E(A), \\
&f \in \mathcal E_1(B), n \in B
\end{aligned} \\
&\begin{aligned}
\beta_{e, f, n} \geq y_{e,f} + q_{n,f} - 1,  \quad \forall &(A, B) \in \mathcal Z, e \in \mathcal E(A), \\
&f \in \mathcal E_1(B), n \in B \label{mccormick-last}
\end{aligned}
\end{align}
and replace the $y_{e,f}q_{n,f}$ terms in (\ref{inter-layer-communication-1}) and (\ref{inter-layer-communication-3}) with $\beta_{e,f,n}$

For fixed network design decisions $x$ and $y$, $\mathcal A(x,y)$ is defined by:
\begin{align}
&\sum_{e \in \mathcal E} z_e \leq U \label{attack-budget}\\
&z_f \leq \sum_{e \in \mathcal E(A)} y_{e,f}z_e, \qquad \forall (A, B) \in \mathcal Z, \forall f \in \mathcal E(B) \label{attack-sequence}\\
& \delta_r = \sum_{e \in \mathcal E(\mathcal S)} x_{e,r} z_e, \qquad \forall r \in \mathcal R \label{relay-attacked}\\
& v_k \leq (1 - \delta_r), \qquad \forall k \in \mathcal K, r \in \mathcal R_k \label{line-off}\\
& v_k \geq \sum_{r \in \mathcal R_k} ( 1 - \delta_r) - |\mathcal R_k| + 1, \qquad \forall k \in \mathcal K \label{line-on}\\
& w_g \leq (1 - \delta_r), \qquad \forall g \in \mathcal G, r \in \mathcal R_g \label{gen-off}\\
& w_g \geq \sum_{r \in \mathcal R_g} ( 1 - \delta_r) - |\mathcal R_g| + 1, \qquad \forall g \in \mathcal G \label{gen-on}\\
& u_d \leq (1 - \delta_r), \qquad \forall d \in \mathcal L, r \in \mathcal R_d \label{load-off}\\
& u_d \geq \sum_{r \in \mathcal R_d} ( 1 - \delta_r) - |\mathcal R_d| + 1, \qquad \forall d \in \mathcal L \label{load-on}
\end{align}
Constraint (\ref{attack-budget}) enforces the attacker's budget. Constraint (\ref{attack-sequence}) requires that the attacker can only access enclaves which are controlled by already-accessed enclaves.
That is, with respect to the graphs in Fig.~\ref{fig:9bus-noattack}, the attack must begin at a balancing authority enclave, and no child can be accessed if its parent was not.
Constraint (\ref{relay-attacked}) enforces that a relay is compromised when the enclave that controls it has been accessed by the attacker. Constraints (\ref{attack-sequence}) and (\ref{relay-attacked}) contain bilinear
terms that can be linearized using the same technique that linearizes (\ref{inter-layer-communication-1}) and (\ref{inter-layer-communication-3}).

Constraints (\ref{line-off})-(\ref{line-on}) enforce that the line $k$ is opened if and only if some relay that controls it is compromised.
Constraints (\ref{gen-off})-(\ref{load-on}) do the same for generators and loads.

Finally, given the attack decision, we define $\mathcal O(u, v, w)$, the operator's DCOPF feasible region:
\begin{align}
&\begin{aligned}
&\sum_{k \in \{k'|d(k') = s\}} f_k - \sum_{k \in \{k'|o(k') = s\}} f_k \\
&\qquad  + \sum_{g \in \mathcal G_s} p_g = \sum_{d \in \mathcal{L}_s} (D_d - l_d)
\end{aligned} & \forall s \in \mathcal S \label{balance}\\
&f_k = B_kv_k(\theta_{o(k)} - \theta_{d(k)} - \Theta_k) & \forall k \in \mathcal K \label{ohms-law}\\
&  -\overline F_k \leq f_k \leq \overline F_k & \forall k \in \mathcal K \label{flow-bounds}\\
& 0 \leq p_g \leq w_g \overline P_g & \forall g \in \mathcal G \label{generation-bounds} \\
&(1 - u_d) D_d \leq l_d \leq D_d & \forall d \in \mathcal L \label{load-shed-bounds} \\
& -\pi \leq \theta_s \leq \pi & \forall s \in \mathcal S \label{angle-bounds}
\end{align}
Equation (\ref{balance}) enforces flow balance at each substation, while the line power flow approximation is enforced in (\ref{ohms-law}).
We linearize (\ref{ohms-law}) by replacing it with:
\begin{align}
&\begin{aligned}
f_k \leq &B_k(\theta_{o(k)} - \theta_{d(k)} - \Theta_k) + \\
&B_k(2\pi + \Theta_k)(1 - v_k)
\end{aligned} &\forall k \in \mathcal K \label{ohms-linearized1}\\
&\begin{aligned}
f_k \geq &B_k(\theta_{o(k)} - \theta_{d(k)} - \Theta_k) \\
&- B_k(2\pi + \Theta_k)(1 - v_k)
\end{aligned} &\forall k \in \mathcal K \label{ohms-linearized2} \\
&- \overline F_kv_k \leq f_k \leq \overline F_kv_k &\forall k \in \mathcal K \label{flow-bounds-linearized}
\end{align}
Constraints (\ref{flow-bounds}) and (\ref{generation-bounds}) enforce thermal limits and maximum generation capacity, forcing flow or generation to 0 if the component is compromised.
Note that we can drop (\ref{flow-bounds}) after adding (\ref{flow-bounds-linearized}).
Also note that in (\ref{generation-bounds}), we assume the minimum generation capacity for all generators is 0.
While this is not accurate, it is necessary in order to make (\ref{balance})-(\ref{angle-bounds}) feasible for all possible values of $u$, $v$, and $w$, the importance of which we describe in Section \ref{sec:methodology}.
In (\ref{load-shed-bounds}), we upper bound load shed by the total load and require that we shed all compromised loads.
Bounds on the phase angles are enforced in (\ref{angle-bounds}).

\section{Solution Methodology} \label{sec:methodology}
Because the grid operator's problem is linear and has an objective opposite the attacker's objective, we can take its dual and linearize the bilinearities that appear in the dual objective function in order to reformulate the inner two levels as one mixed integer linear program (MILP).
For details on this procedure, see \cite{motto2005mixed}.
Note that the linearization requires upper bounds on the dual variables.
The authors of \cite{KleinertLPS2019} show that verifying the correctness of these dual bounds is as hard as solving the original bilevel problem.
Thus, we adopt the heuristic commonly used in prior literature, which is to assume that the duals are bounded above by the power capacity of the component corresponding to their index \cite{salmeron2009worst}. 
After this reformulation, where we take the dual of the linear program which minimizes $\sum_{d \in \mathcal L} l_d$ over the set $\{(\theta, f, p, l) : \text{(\ref{balance}), (\ref{generation-bounds})-(\ref{angle-bounds}), (\ref{ohms-linearized1})-(\ref{flow-bounds-linearized})}\}$ and combine the resulting maximization problem with the attacker problem,
we have a bilevel problem with integer leader and mixed integer follower:
\begin{align}
\min_{(x, y, q, t, \beta) \in \mathcal D} \quad \max_{\bm X \in \mathcal A_D(x,y)} L, \label{bilevel-reformulation}
\end{align}
where 
\begin{align*}\bm X = (&\delta, z, u, v, w, f, p, d, \theta, 
\mu, \xi^+, \xi^-,  \\ &\lambda^+,\lambda^-, \gamma, \alpha^+, \alpha^-, \beta^+, \beta^-)
\end{align*}
and $\mathcal A_D(x, y)$ is defined by:
\begin{align}
&\text{(\ref{attack-budget})-(\ref{load-on})} \nonumber \\
&\begin{aligned}
&\sum_{s \in \mathcal S} \big [\mu_s \sum_{d \in \mathcal{L}_s} D_d - \pi(\beta^+_s + \beta^-_s) \big ] \\
&+ \sum_{d \in \mathcal L} D_d[\alpha^+_d(1 - u_d) - \alpha^-_d]  \\
&- \sum_{k \in \mathcal K} [B_k(\Theta_k + (2\pi + \Theta_k)(1 - v_k))(\xi^+_k + \xi^-_k) \\
& \qquad \quad + \overline F_kv_k (\lambda^+_k + \lambda^-_k)] - \sum_{g \in \mathcal G} \overline P_g \gamma_g w_g = L \\
\end{aligned} \label{strong-duality} \\
&\mu_{d(k)} - \mu_{o(k)} + \xi^+_k - \xi^-_k + \lambda^+_k - \lambda^-_k = 0, \quad \forall k \in \mathcal K \label{flow-dual}\\
&\mu_{s(g)} - \gamma_g \leq 0, \quad \forall g \in \mathcal G \label{generation-dual}\\
&\mu_{s(d)} + \alpha^+_d - \alpha^-_d \leq 1, \quad \forall d \in \mathcal L \label{load-shed-dual}\\
&\begin{aligned}
\beta^+_s - \beta^-_s &+ \sum_{k \in \{k'| o(k') = s\}} B_k(\xi^-_k - \xi^+_k) \\
& + \sum_{k \in \{k'| d(k') = s\}} B_k(\xi^+_k - \xi^-_k)  = 0,
\end{aligned} \quad \forall s \in \mathcal S \label{angle-dual} \\
&\xi^+, \xi^-, \lambda^+,\lambda^-, \gamma, \alpha^+, \alpha^-, \beta^+, \beta^- \geq 0 \label{dual-bounds}
\end{align}
where $\mu$ is the dual of the balance constraint (\ref{balance}), $\xi^-$ and $\xi^+$ are the duals of (\ref{ohms-linearized1}) and (\ref{ohms-linearized2}) respectively, $\lambda^+$ and $\lambda^-$ are the duals of (\ref{flow-bounds}), $\gamma$ is the dual of the upper bound in (\ref{generation-bounds}), $\alpha^+$ and $\alpha^-$ are the duals of (\ref{load-shed-bounds}), and $\beta^+$ and $\beta^-$ are the duals of (\ref{angle-bounds}).
Constraint (\ref{strong-duality}) calculates the dual objective value.
Constraints (\ref{flow-dual}), (\ref{generation-dual}), (\ref{load-shed-dual}), and (\ref{angle-dual}) are the dual constraints corresponding to $f$, $p$, $l$, and $\theta$, respectively.
Constraint (\ref{strong-duality}) includes many bilinear terms, all the product of a non-negative dual variable and a binary variable.
These can be reformulated with the addition of auxiliary continuous variables.
For example, to reformulate $\alpha^+_d u_d$, we introduce $\overline \alpha^+_d$ in place of the product in constraint (\ref{strong-duality}) and add the constraints
\begin{align}
&\alpha^+_d - D_d(1 - u_d) \leq \overline \alpha^+_d \leq \alpha^+_d + D_d(1 - u_d) \\
&0 \leq \overline \alpha^+_d \leq D_d u_d,
\end{align}
where we are assuming that $\alpha^+_d$ is bounded above by the total demand at load $d$.

Note that for every outer-problem solution $(x, y, q, t, \beta) \in \mathcal D$, problem (\ref{bilevel-reformulation}) has a finite objective value.
That is, $\mathcal A_D(x,y) \ne \emptyset$.
This is because, regardless of the network segmentation and the attack, it is always possible to prevent infeasibility by shedding load or turning off generators.
We refer to this property as {\it relatively complete recourse.}
Note that it is because the generator dispatch lower bounds are 0 that we have this property.

Problem (\ref{bilevel-reformulation}) is a mixed-integer bilevel model, which we solve using the bilevel branch and bound algorithm from \cite{Fischetti2018}. 
In \cite{Fischetti2016}, \cite{Fischetti2017}, and \cite{Fischetti2018}, callbacks within the IBM CPLEX solver \cite{cplex} are used to create a bilevel branch-and-cut solver. 
The authors make their software available for academic use at \cite{fischetti_solver}. 
To easily generate the required MPS file for the bilevel branch and bound solver, we formulated our model with Pyomo \cite{pyomo-book} and \cite{pyomo-paper}.

\subsection{An Equivalent Solution Methodology}\label{sec:alt-methodology}
We found that, when solving (\ref{bilevel-reformulation}) with the solver from \cite{fischetti_solver}, the solver could not close the gap within a time limit of 72 hours. 
We suspect this may be a symptom of numerical issues in the cuts.
To overcome this, we solved a slightly modified problem, in which the variable $\beta$ and the constraints it appears in are moved to the inner problem of (\ref{bilevel-reformulation}).
Formally, the formulation is as follows:
\begin{equation}\label{hack-bilevel}
\min_{(x,y,q,t) \in \mathcal D'} \max_{(\bm X, \beta) \in \mathcal A_D'(x,y,q)} L
\end{equation}
where $\mathcal D' = \{(x, y, q, t): \text{(\ref{substation-subnets-talk-to-a-relay})-(\ref{each-subnet-assigned-one-entity}), (\ref{inter-layer-communication-2}), (\ref{communication-implies-belonging-existing})-(\ref{t-domain})}\}$ and
\[
\mathcal A_D'(x, y, q, t) = \{(\bm X, \beta) : \text{(\ref{mccormick1})-(\ref{load-on}), (\ref{strong-duality})-(\ref{dual-bounds})} \}.
\]
Note that, in (\ref{hack-bilevel}), we have lost the relatively complete recourse property: Fixing the outer level problem's choice of $y$ and $q$, it is possible to choose $\beta$ in order to make one or more of (\ref{mccormick1})-(\ref{mccormick-last}) infeasible.
However, because the solver from \cite{fischetti_solver} uses the high-point relaxation of (\ref{hack-bilevel}) to get a lower bound, it assumes the outer problem will not choose to make the inner problem infeasible.
That is, we are actually solving the problem:
\begin{equation}\label{hpr-bilevel}
\min_{(x,y,q,t) \in \mathcal P} \max_{(\bm X, \beta) \in \mathcal A_D'(x,y,q)} L,
\end{equation}
where
\[
\mathcal P = \{(x,y,q,t) \in \mathcal D': \exists \, (\bm X, \beta) \text{ s.t. } (\bm X, \beta) \in \mathcal A_D'(x,y, q)\}
\]
Essentially, there is an implicit constraint added to the outer problem that its solution is in the projection of the inner problem's feasible region onto the outer problem variables.
We therefore show that (\ref{hpr-bilevel}) is equivalent to (\ref{bilevel-reformulation}).

We first show that the optimal objective value of (\ref{hpr-bilevel}) is a lower bound to that of (\ref{bilevel-reformulation}).
Note first that $\mathcal P$ equals the projection of $\mathcal D$ onto the $(x, y, q, t)$ variables.
This follows from the relatively complete recourse property of (\ref{bilevel-reformulation}).
 We can rewrite $\mathcal P$ as
\[
\mathcal P = \{(x, y, q, t) \in \mathcal D' : \exists \, \beta \text{ s.t. (\ref{mccormick1})-(\ref{mccormick-last})}\}
\]
Next, note that the projection of $\mathcal A_D'(x, y, q)$ onto the space of the $\bm X$ variables is a restriction of $\mathcal A_D(x,y)$ (since constraints (\ref{mccormick1})-(\ref{mccormick-last}) have been moved to the inner problem).
Since, in projected space, we have an equivalent outer problem region and a restricted inner problem region, we have a lower bound.
Next, we show that this bound is achieved.
To see this, consider a solution $(x, y, q, t, \beta, \bm X)$ to (\ref{hpr-bilevel}).
Then, by our argument above, $(x, y, q, t, \beta) \in \mathcal D$, and by the definition of $\mathcal A_D'(x, y, q)$, $\bm X \in \mathcal A_D(x,y)$.
Thus, this solution is feasible in (\ref{bilevel-reformulation}).
This completes the argument: 
When using the branch and bound solver from \cite{fischetti_solver} to solve (\ref{hack-bilevel}), we are solving (\ref{bilevel-reformulation}) exactly.

\section{Results}\label{sec:results}

We present a case study on the 9-bus WSCC and 30-bus IEEE systems \cite{9_bus} using the methodology given in Section \ref{sec:alt-methodology}.
Note that, as we have formulated the network segmentation problem, the network designer is forced to use their entire budget, regardless of if every enclave in the budget is necessary to reduce that load shed the attacker can cause.
While it would be possible to reformulate the problem to avoid this, it would be necessary to use the network designer's objective to motivate them to be frugal with their budget.
This is no longer an interdiction problem, and is in essence a multi-objective variant.
For simplicity, we leave this as future work, and instead suggest solving the proposed model for a variety of network designer budgets in order to find a minimal one which achieves a satisfactory reduction in load shed.
We will show an example of this process for the 9-bus system.

\subsection{Hardware and Software Specification}
We ran our experiments on a Linux server with two Intel\textregistered \; Xeon\textregistered \; Silver 4210 CPU @ 2.20GHz processors and 256 GB RAM. 
For the CPLEX-based solver from \cite{fischetti_solver}, we used CPLEX 12.9.0 with two threads.

\subsection{9-bus WSCC System}\label{subsec:9-bus-results}
\begin{figure*}
	\centering	
	\subfloat[The 9-bus system before segmentation. The worst-case attack with budget 5 is marked with a skull-and-crossbones on each enclave and component which is compromised. The attacker infiltrates the Balancing Authority enclave, the Control Center 2 enclave, and the three substation enclaves whose relays control the loads and lines at Buses 5, 6, and 8. ]{\includegraphics[width=0.75\linewidth]{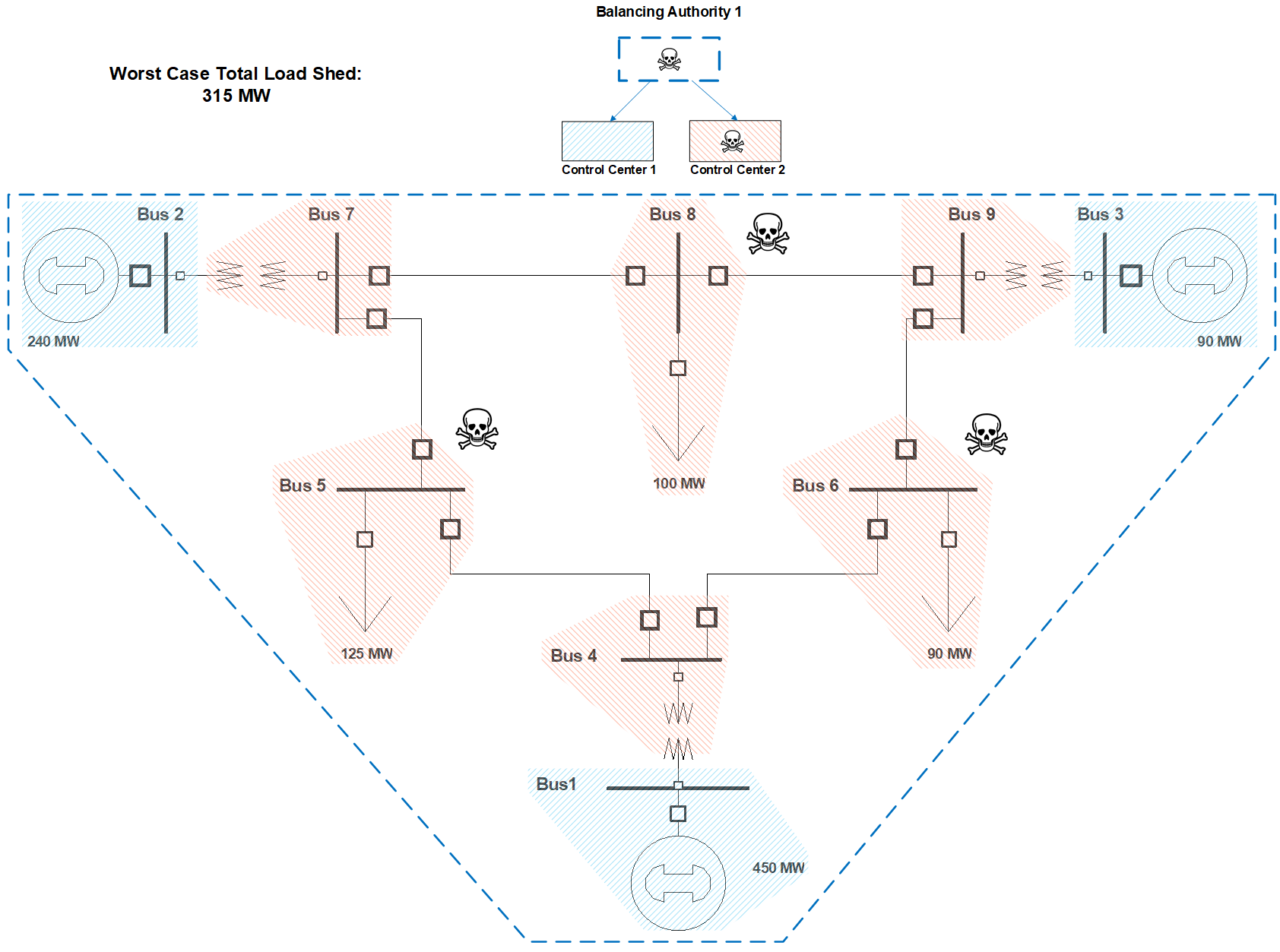}\label{fig:9bus-network-before-segmentation-budget5}}\\
	\subfloat[The 9-bus system after segmentation, where the network designer was allowed a budget of two control center enclaves and the attacker had a budget of five enclaves.
	After segmentation, the attacker spends three units of budgets to infiltrate the balancing authority  enclave and both the Control Center 2 enclaves. They only have two units of budget remaining, so they infiltrate buses 4 and 7, which islands the two largest generators, resulting in 225 MW of load shed. Note that the attack on the second yellow Control Center 2 enclave does not cause any load shed; it only allows the attacker to use the full attack budget.]{\includegraphics[width=0.75\linewidth]{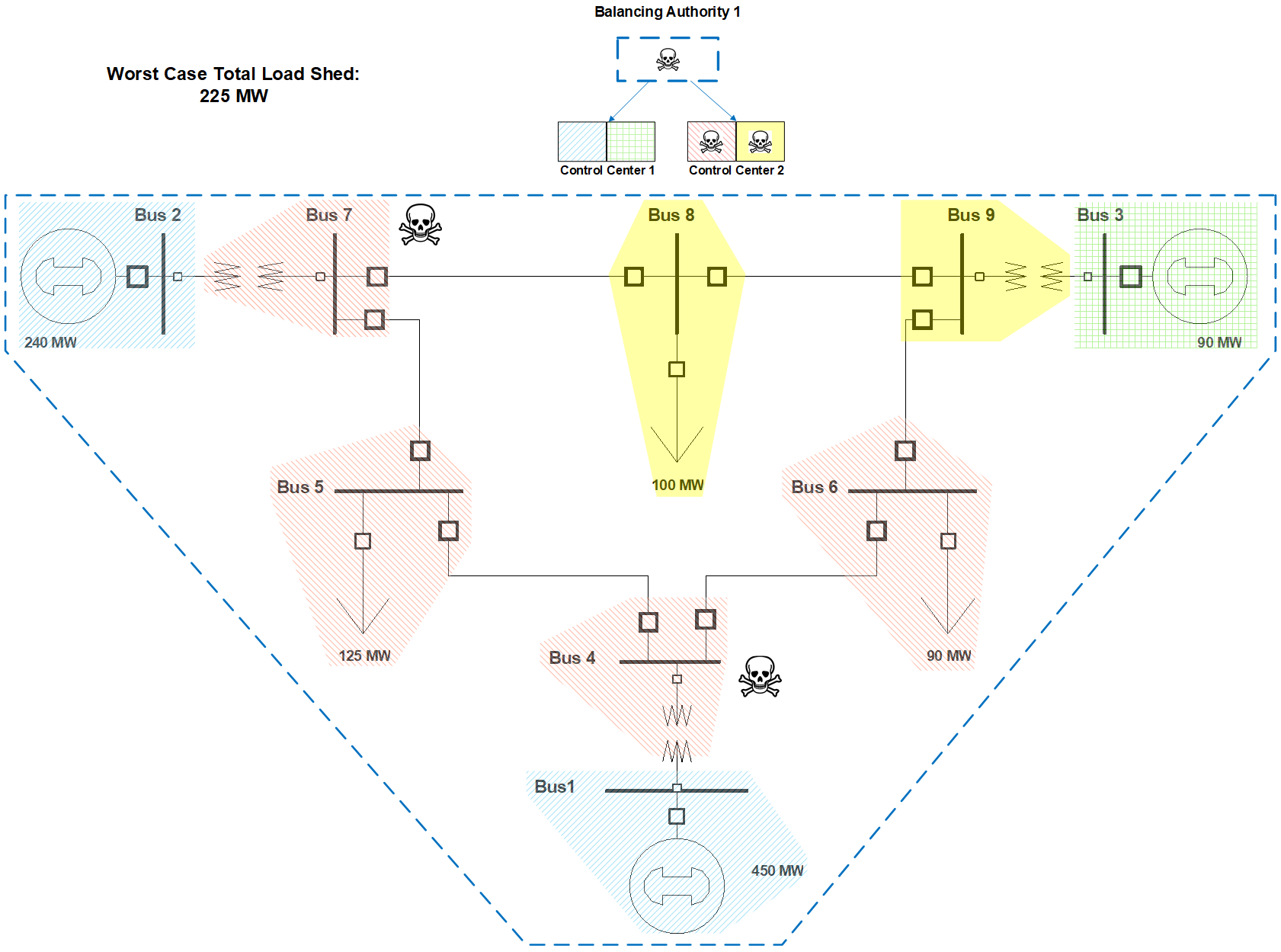}\label{fig:9bus-network-after-segmentation-budget5}}
	\caption{The 9-bus system before and after segmentation in figures (a) and (b) respectively, when the attacker has a budget of 5 enclaves.}
	\label{fig:9bus}
\end{figure*}

To create intuitive results that are easy to validate, we first use the 3-generator 9-bus WSCC test system \cite{9_bus}.
We use a simple communication network consisting of a single balancing authority and two control centers. 
The three substations with generators are assigned to one control center and the remaining six substations are assigned to the other control center. 
There is exactly one security enclave at each substation and each grid component at that substation is controlled through a relay hosted by its substation's security enclave. 
See Fig.~\ref{fig:9bus-network-before-segmentation-budget5} for a depiction of this cyber-physical system before segmentation.
Note that this starting topology does not correspond exactly to the communication network we showed in Fig.~\ref{fig:9bus-noattack}.
In Fig.~\ref{fig:9bus-noattack}, we reduced the number of relays in the system for the sake of readability.

\begin{table}
	\caption{Partial paramter sweep of network designer budgets for the 9-bus system with attack budget 5.}\label{table:9bus-budget5-parameter-sweep}
\begin{tabular}{r | r | r | r}
	\hline
	\makecell{No. of \\ Substation \\ Enclaves} & \makecell{No. of \\ Control Center\\ Enclaves} & \makecell{No. of Balancing \\ Authority \\ Enclaves} & \makecell{Worst-Case Attack \\ Load Shed (MW)} \\
	\hline
	0 & 1 & 0 & 315 \\
	1 & 0 & 0 & 315 \\
	1 & 1 & 1 & 315\\
	0 & 2 & 0 & 225 \\
	0 & 0 & 2 & 190 \\
	\hline
\end{tabular}
\end{table}

We do two case studies on this network: The first defends against an attacker with a budget of 5 and the second defends against an attacker with a budget of 8.
We begin with the case where the budget is 5. Before network segmentation, the attacker can cause a complete blackout by pivoting from the balancing authority, through the second control center, to Substations 5, 6, and 8. 
These substations contain all of the system's load of 315 MW. 
See Fig.~\ref{fig:9bus-network-before-segmentation-budget5} for a visualization of this attack. 
Note that the attacker has multiple optimal solutions: For example, they could alternatively pivot from the balancing authority to Control Center 1, and then compromise Substations 1, 2, and 3.
This allows them to shut down all the generators, also resulting in 315 MW of load shed.

In Table~\ref{table:9bus-budget5-parameter-sweep}, we present results on a partial parameter sweep of network designer budgets for this problem.
We see that two enclaves at either the balancing authority or control center level are necessary to prevent any load shed.
Furthermore, two enclaves at the balancing authority level prevent 35 MW more load shed than two enclaves at the control center level. 

Figure \ref{fig:9bus-network-after-segmentation-budget5} displays the optimal network segmentation for a budget of two control center enclaves.
The model segments both the control center networks into two enclaves, and reallocates the substations so that it is no longer possible to compromise all three generators or all three loads. 
After network segmentation, one worst-case attack infiltrates both enclaves in Control Center 2. 
The attacker can then only access two substations.
They choose Substations 4 and 7 to cause a load shed of 225 MW, reducing the load shed possible in the unsegmented network by 28.6\%. 
Note that the attack on the solid yellow right enclave of Control Center 2 is unnecessary.
However, due to the segmentation, there is nothing productive for the attacker to do with this last unit of budget.

\begin{figure}
	\centering
	\includegraphics[width=\linewidth]{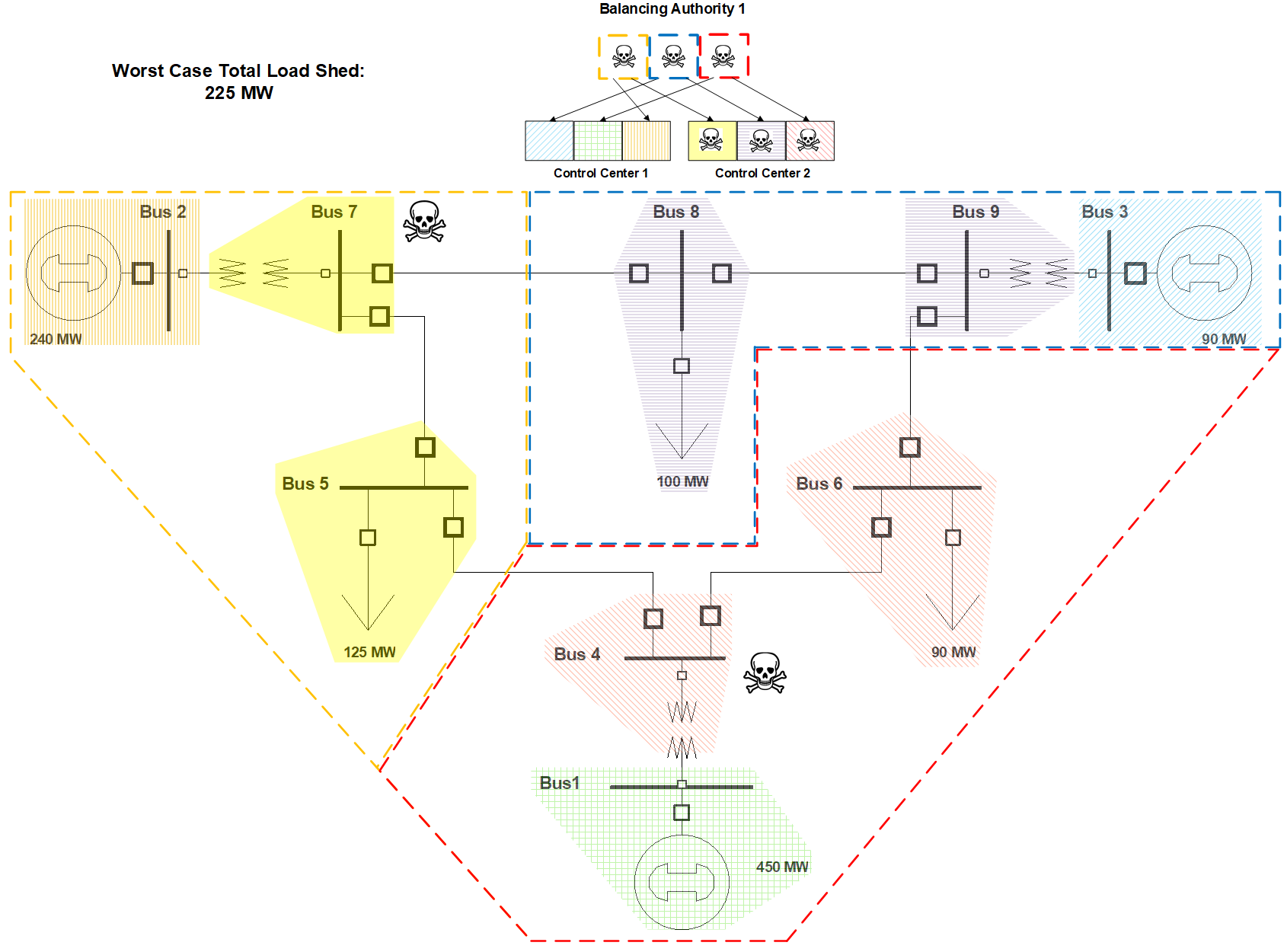}
	\caption{The 9-bus system after segmentation with a network designer budget of 2 balancing authority enclaves and 4 control center enclaves and an attack budget of 8 enclaves.}\label{fig:9bus-network-after-segmentation-budget8}
\end{figure}

In the second case study, we decide an optimal network segmentation to protect against an attacker with a budget of 8 enclaves.
Clearly, since the attacker with budget 5 was able to shed all the load (see Figure~\ref{fig:9bus-network-before-segmentation-budget5}), in the unsegmented case with budget 8, the attacker sheds 315 MW of load also.
Using a process similar to what we showed in Table~\ref{table:9bus-budget5-parameter-sweep}, we select a network designer budget of two balancing authority enclaves and four control center enclaves.
The resulting segmentation is visualized in Fig.~\ref{fig:9bus-network-after-segmentation-budget8}.
The balancing authority is divided into three enclaves, each of which communicates with two control center enclaves.
No load bus or generator bus is on the same control center enclave, so, to remove all the loads or all the generators, the attacker would need to spend three units of budget at the substation level and three at the control center level.
However, this leaves them only two at the balancing authority level.
Because each balancing authority enclave is mapped to one Control Center 1 enclave and one Control Center 2 enclave, they would need nine or more units of budget to attack enough enclaves to cause a complete blackout.
Thus, they cannot shed all the load, and 90 MW (the smallest load) is served.
Again, network segmentation reduced the load shed by 28.6\%.
Note in the attack shown in Fig.~\ref{fig:9bus-network-after-segmentation-budget8}, the attacker wastes two units of budget: Similarly to the budget-5 case, the attacks on the blue balancing authority enclave and the purple control center enclaves are unproductive, but they cannot be reassigned so that the attack increases the load shed.

Both of these examples illustrate that optimal network segmentation isolates the most critical infrastructure in a way that drastically reduces the damage caused by the worst-case attack.

\subsection{30-bus System}

\begin{figure*}
	\centering
	\subfloat[The 30-bus system before segmentation: The network is controlled by two balancing authorities and 3 control centers, each with only one security enclave. The worst-case attack with budget 6 compromises Balancing Authority 1, Control Center 1, and four substation enclaves controlling buses Blaine, Claytor, Kumis, and Reusens.]{\includegraphics[width=0.85\linewidth]{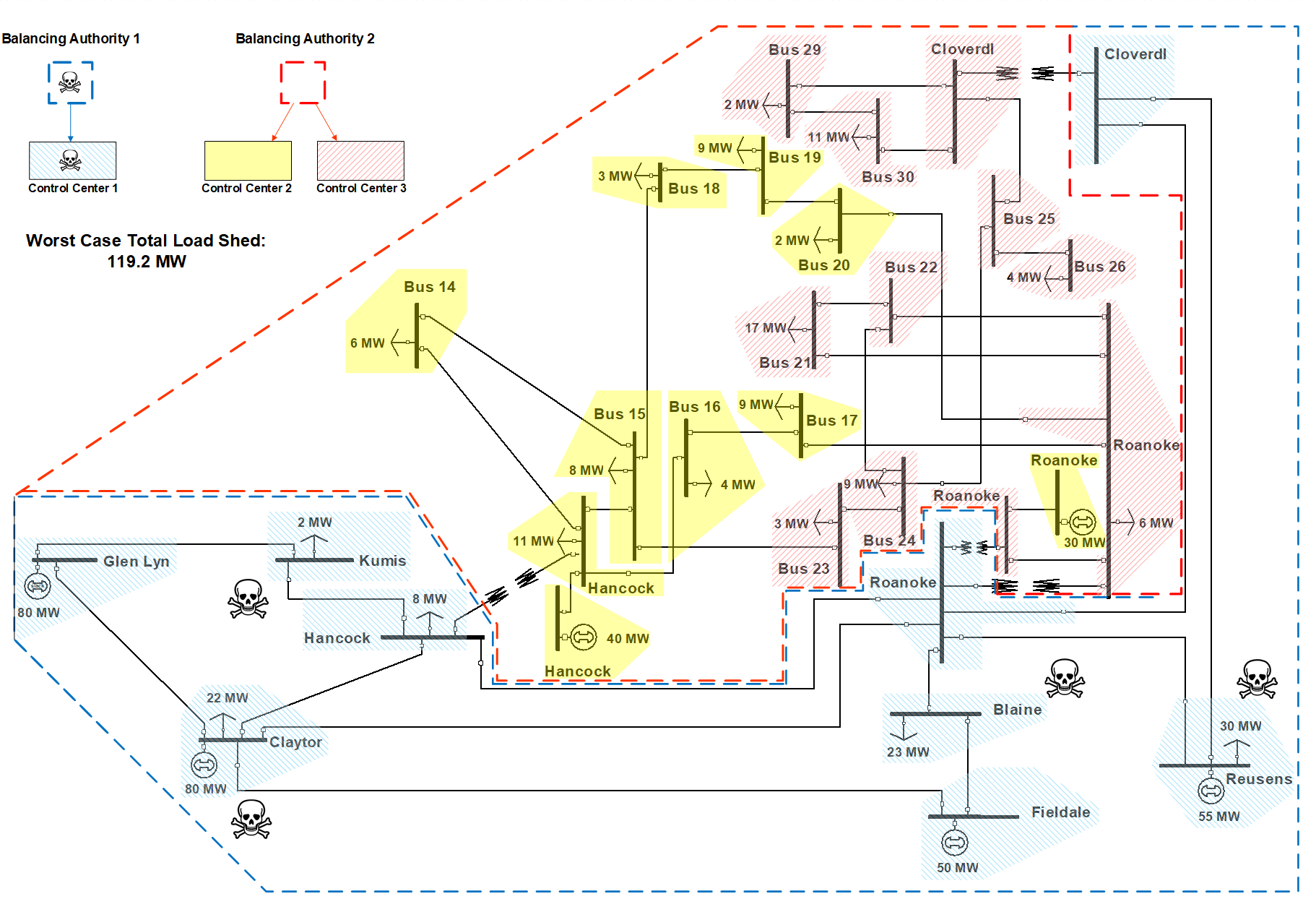}\label{fig:30bus-before-segmentation}} \\
	\subfloat[The 30-bus system after segmentation: Balancing Authority 1 and Control Center 1 have two enclaves each, meaning the attacker would require 4 units of budget to gain access to the same substations. In lieu of this, they attack only the green control center enclave, compromising Blaine, Fieldale, and Reusens. While it is possible for them to also attack Glen Lyn or Cloverdl, the resulting load shed is not higher.]{\includegraphics[width=0.85\linewidth]{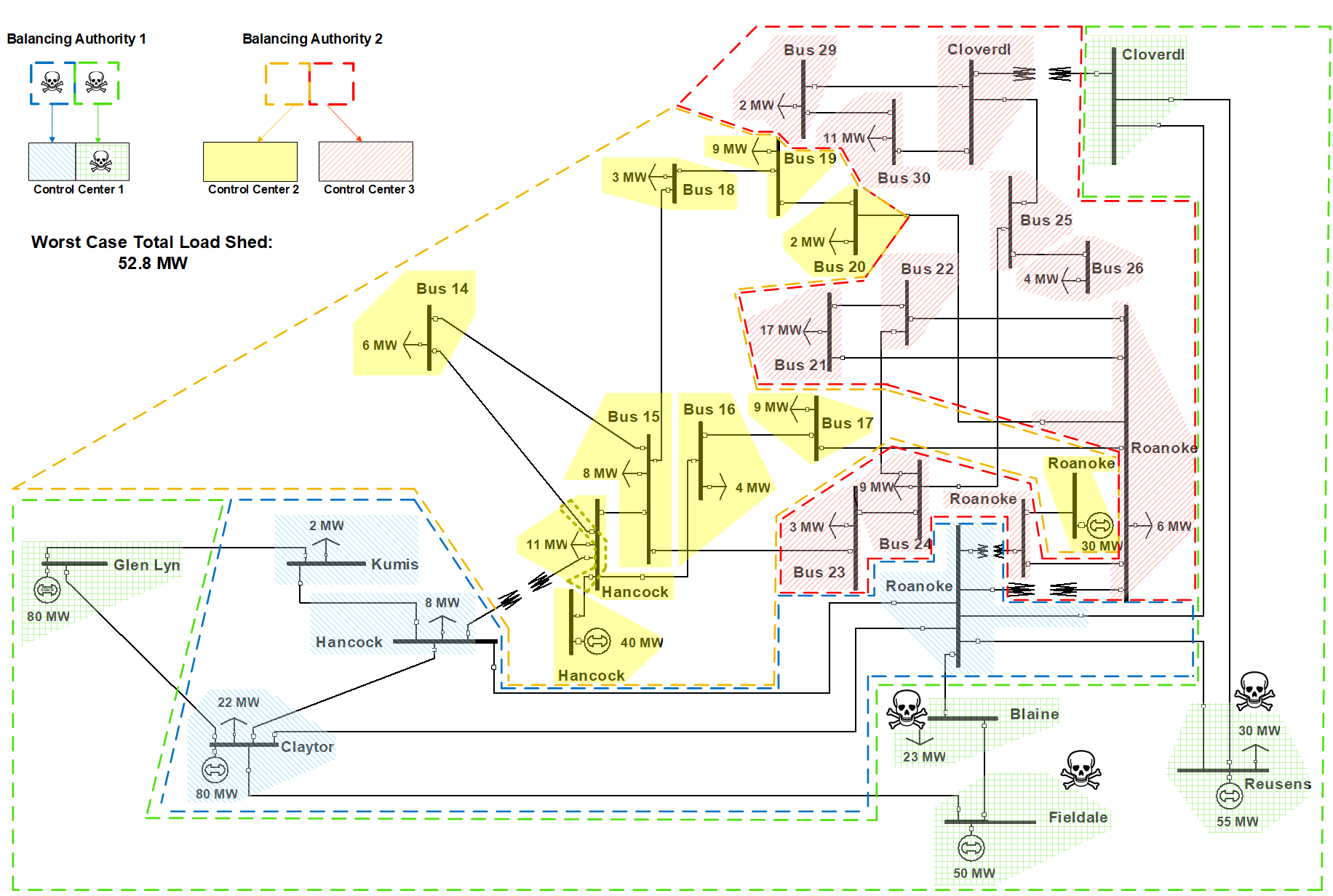}\label{fig:30bus-after-segmentation}}
	\caption{The 30-bus system before and after segmentation in figures (a) and (b) respectively, with the worst-case attack with budget 6 marked in both versions.}
	\label{fig:30bus-segmentation-results}
\end{figure*}

The 30-bus IEEE system provides a larger example with a more complex topology.
For this system, we use the voltage levels to create a fictitious communication network. 
All high-voltage grid components communicate with a single control center and balancing authority. 
The lower-voltage grid components are split between two control centers.
Both these control centers communicate with a second balancing authority.
This network is depicted in Fig.~\ref{fig:30bus-before-segmentation}.
Balancing Authority 1 and Control Center 1 control the high-voltage grid components, and the remaining two control centers and Balancing Authority 2 control the low-voltage components.

In this case study, we consider an attack budget of six enclaves.
We again search the results for a range of network designer budgets, trying to find a minimal one with satisfactory reduction of load shed.
We present results where we allow two extra balancing authority enclaves, one extra control center enclave, and one extra substation enclave.
Our model segments both balancing authorities into two security enclaves each and divides Control Center 1 into two enclaves. 
Because substations controlled by Control Center 1 have the largest loads and generation, the model assigns substations to the two enclaves in Control Center 1 so that the largest loads communicate with only one control center enclave.
Control Centers 2 and 3 are each assigned to their own enclave in Balancing Authority 2.
Finally, the substation enclave is used at the Hancock bus in Control Center 2 with 11 MW of load: The relays controlling the lines to Bus 14 and to the Hancock bus with the 40 MW generator are on their own enclave.

Before network segmentation, the attacker can cause a load shed of 119.2 MW, or 63.1\% of the total load, by pivoting from the first balancing authority, through its control center, to substations Kumis, Claytor, Blaine, and Reusens. See Fig.~\ref{fig:30bus-before-segmentation}. 
After network segmentation, the worst-case attack uses two units of budget to give the attacker access to the substations controlled by the green control center enclave. 
This allows them to attack substations Blaine, Fieldale, and Reusens, causing a total load shed of 52.8 MW, or 27.9\% of the total demand.
This means that network segmentation reduced the worst-case load shed by 55.7\%. 
See Fig \ref{fig:30bus-after-segmentation} for a depiction.
Note that the attack on the blue balancing authority enclave is wasted.
In fact, the attack on Fieldale is also not necessary: All of the load shed stems from directly compromising the loads at Blaine and Reusens, and the generators at Glen Lyn and Fieldale are not necessary for serving the remaining load.

This case study illustrates an intuitive network segmentation strategy. The largest loads are connected to substations Blaine and Reusens, so an attacker will expend all resources to infiltrate those substations. 
However, this network segmentation strategy ensures that the attacker will not be able to shed any further load.
While it may seem intuitive to put Blaine and Reusens on separate control center enclaves, note that the attacker still has enough budget to access both, expending four units of budget to attack both Balancing Authority 1 enclaves and both Control Center 1 enclaves, then using the remaining two units to attack Blaine and Reusens. 
The segmentation therefore limits the attacker's ability to compromise additional enclaves.

Finally, we note that the segmentation of Balancing Authority 2 and of the Hancock substation are indeed necessary in this example.
Without segmenting Balancing Authority 2, there is an attack on 6 enclaves which sheds 104.7 MW of load by compromising the Balancing Authority 2 enclave, the enclaves in Control Centers 1 and 2, the Roanoke substation with the 6 MW load, and the Cloverdl substation under Control Center 3.
By segmenting Balancing Authority 2, this attack is prevented since it takes an additional unit of budget to gain access to both control centers.
With only this segmentation and not the extra substation enclave at Hancock, there still exists an attack which shed 53 MW of load, targeting the yellow Balancing Authority 2 enclave, Control Center 2, the Hancock substation with the 11 MW load, Bus 15, Bus 17, and Bus 20.
Including the new Hancock substation enclave, the worst attack on the low-voltage side of the network is reduced to 51.7 MW, compromising the red Balancing Authority 2 enclave, Control Center 3, the Roanoke substation with the 6~MW load, Bus 23, and the accessible Cloverdl substation.

This case study shows that even moderate budgets of additional enclaves can greatly decrease the possible load shed achievable by a cyber attacker, in this case reducing it by more than half.

\section{Future Work}\label{sec:future-work}
The network segmentation model we develop can be improved and extended in several ways. 
More details from the communication network presented in \cite{gaudet2020firewall} can be included to result in a more realistic network. 
The attacker model can be made more realistic by adding more complex movement throughout the network, using enclave access fees based on real data, and adding context to the attacker's budget. For example, an IT administrator could perform analysis to determine the maximum number of enclaves that could likely be comprised before the IDS systems on the network would generate an alert, and use that as the basis for making segmentation decisions. 
We use the simplest grid operator model available. 
Future work can focus on using a higher-fidelity power flow model and allowing grid operator transmission switching.
In addition, solution strategies that would allow solving the model on realistically-sized networks would be of interest.

\section{Conclusion}\label{sec:conclusion}
Segmentation is a cyber defense strategy that has been proposed for improving network security. However, specific strategies for optimally implementing these recommendations, while accounting for the underlying system being protected, are lacking. To address this, we develop a trilevel cyber-physical power-transmission-system network-segmentation model that is the first of its kind. It models an IT administrator who must decide how to preemptively segment the power transmission system communication network in preparation for an attacker with perfect information of the whole system. 
We use an attacker model where a malicious actor must move from balancing authorities through control centers to substations in order to disable grid components through relays. 
The attacker anticipates how the grid operator will use a DCOPF to redispatch the generators after the attack is executed. 
We reformulate our model to a mixed-integer bilevel problem and use bilevel branch-and-cut to solve it. 
Finally, we present network segmentation results on the 9-bus WSCC and 30-bus IEEE test systems, demonstrating the benefit of optimal network segmentation in terms of reducing the worst-case load shed the attacker can achieve.

\section{Acknowledgments}

We would like to thank Santanu Dey for his help in showing the correctness of the methodology in Section \ref{sec:alt-methodology}.

\bibliographystyle{IEEEtran}
\bibliography{network_segmentation_paper}{}

\end{document}